\documentclass[a4paper,10pt]{article}
\usepackage[utf8]{inputenc}
\usepackage{amssymb,graphicx,epsfig,amsmath,amsfonts}
\usepackage{color}
\usepackage{colortbl}
\usepackage{algorithm}
\usepackage{algpseudocode}

\usepackage{amssymb,graphicx,epsfig,amsmath,amsfonts}
\usepackage{color,algorithm,algpseudocode}
\newcommand{\Nn}{{\mathbb{N}}}

\def\a2{{\alpha/2}}

\def\ga{{\gamma}}
\def\be{{\beta}}

\def\bfq{{\mathbf{q}}}

\def\bfz{{\mathbf{z}}}

\def\a2{{\alpha/2}}

\def\GG{{\Gamma}}

\def\bfzeta{\mbox{\boldmath$\zeta$}}

\def\qed{\hfill$\Box$}
\def\cR{{\cal{R}}}
\def\cP{{\cal{P}}}

\newtheorem{theorem}{Theorem}[section]
\newtheorem{remark}{Remark}[section]
\newtheorem{proposition}{Proposition}[section]

\title{Weakly regular Sturm-Liouville problems: a corrected spectral matrix method.}

\author{Cecilia Magherini\thanks{Dipartimento di Matematica, Universit\`{a} di Pisa, Italy, cecilia.magherini@unipi.it.}}

\date{}

\begin{document}

\maketitle

\begin{abstract}
In this paper, we consider weakly regular Sturm-Liouville eigenproblems
with unbounded potential at both endpoints of the domain. 
We propose a Galerkin spectral matrix method for its solution and 
we study the error in the  eigenvalue approximations it provides. 
The result of the convergence analysis is then used 
to derive a low-cost and very effective formula for the computation of corrected numerical eigenvalues.
Finally,  we present and discuss the results of several numerical experiments which  confirm the validity of the approach.

\end{abstract}

\section{Introduction}
Recently, the author studied a corrected spectral matrix method for solving weakly regular and singular Sturm-Liouville 
problems defined over the bounded domain $(-1,1)$ with an unbounded potential at the left endpoint, \cite{Mag}. The numerical 
results provided by such technique are definitely satisfactory  for weakly regular problems. This suggested 
to study a generalization of the method for the approximation of the eigenvalues and of the eigenfunctions 
of problems of the following type 
\begin{eqnarray}\label{slp}
-y''(x) + q(x) y(x) &=& \lambda y(x), \qquad  x \in (-1,1)\>, \\
\alpha_L y(-1) + \beta_L y'(-1) &=& 0\>, \qquad \quad\alpha_L^2+\beta_L^2 \neq 0,\label{bca}\\
\alpha_R y(-1) + \beta_R y'(-1) &=& 0\>, \qquad \quad\alpha_R^2+\beta_R^2 \neq 0,\label{bcb}
\end{eqnarray}
where the potential $q$ is given by
\begin{equation}\label{pot1}
q(x) = \sum_{i=1}^S\frac{g_i(x)}{(1-x)^{\be_i} (1+x)^{\ga_i}}, \qquad \be_i,\ga_i <1, \quad i=1,\ldots,S,
\end{equation}
with functions $g_i$ at the numerators that are analytical inside and on a Bernstein ellipse containing $[-1,1].$ In the literature,
problems of this type with $q$ unbounded at least at one endpoint 
are sometimes called weakly regular and it is well known that their spectrum is composed by real and simple eigenvalues 
which can be ordered as an increasing sequence tending to infinity. We will number them starting from index $k=1,$ i.e. we will call
$$\left\{\lambda_1 < \lambda_2 < \lambda_3< \ldots\right\}$$
the exact spectrum of (\ref{slp})--(\ref{pot1}).\\

In \cite{Mag}, we considered the case $S=2$ with $\beta_1=\gamma_1=\beta_2=0,$  namely problems with a 
potential of the form $q(x) = g_1(x) + g_2(x)/(1+x)^{\ga_2},$  and a special algorithm for $\ga_2\in(0,1)$ and $y(-1) \neq 0$
was derived. 
As remarked in the same paper, the obtained results appear to be more reliable than those given by well-known 
and well-established general-purpose codes based on shooting techniques like, for example, 
the MATSLISE2 \cite{Led,LedTh}, the SLEDGE \cite{Sledge} and the SLEIGN2 \cite{Sleign} ones. 
A possible explanation is that the common basic idea in them implemented is essentially the selection of suitable layers.
In particular, if $q$ is unbounded at both endpoints, then the approach is that of solving a suitable problem over 
$(-1+\varepsilon_1,1-\varepsilon_2)$ with $\varepsilon_1$ and $\varepsilon_2$ small positive values automatically selected, \cite{Pryce}. 
As indicated in their documentation, this may cause a loss of accuracy  
and from all our tests we deduced that this may be more relevant if the problem 
is not subject to the Dirichlet condition at the endpoint where $q$ is unbounded.
As an example, in the following table 
we list some numerical eigenvalues for the problem with $q(x) = 10/(1-x^2)^{3/4}$ 
subject to $y'(\pm 1)=0$ that we computed by using such codes with a tolerance equal to $10^{-13}.$\\

\begin{small}
\begin{center}
\begin{tabular}[h]{|r||r|r|r|}
\hline
\vspace{-.2cm}&&&\\
\multicolumn{1}{|c||}{$k$}&
\multicolumn{1}{|c|}{MATSLISE2} &
\multicolumn{1}{|c|}{SLEDGE} &
\multicolumn{1}{|c|}{SLEIGN2} \\
&&&\vspace{-.25cm}\\
\hline
&&&\vspace{-.25cm}\\
$5$ & $70.95684246381$ & $70.97934056277$ & $70.94410214688$ \\
$10$ & $235.44749698614$ & $235.51215075034$ & $235.40992793209$ \\
$20$ & $925.03586530831$ & $925.11877235794$ & $924.98715263138$ \\
\hline
\end{tabular}
\end{center}
\end{small}

These considerations justify the interest in generalizing the method proposed in \cite{Mag} 
and the outline of this paper is the following. In Section~\ref{sec2}, we recall the basic facts concerning the spectral
Legendre-Galerkin matrix method introduced in \cite{Mag} and we discuss the computation of the coefficient matrix 
that corresponds to a potential $q$ of the form in (\ref{pot1}). 
An analysis of the error in the numerical eigenvalues with respect to the generalized eigenvalue problem size is carried
out in Section~\ref{sec3}. In addition, in the same section, we derive a low cost and effective procedure 
for an a posteriori correction of the numerical eigenvalues. Finally, in Section~\ref{sec4} we report and discuss the 
results of some numerical experiments.
\section{Spectral Legendre-Galerkin method}\label{sec2}
Let $\Pi_{N+1}$ be the space of polynomials of maximum degree $N+1,$ 
for a fixed $N\in \Nn,$ and let
\begin{eqnarray}
{\cal S}_N &\equiv& \left\{  r \in \Pi_{N+1}: \quad 
\alpha_L \,r(-1)+\beta_L \,r'(-1)= \alpha_R \,r(1) + \beta_R\, r(1) = 0 \right\}\label{SN}\\
       &\equiv& \mbox{span} \left(\cR_0,\cR_1,\ldots,\cR_{N-1}\right). \label{spanSN}
             \end{eqnarray}
We look for an approximation of an eigenfunction $y$ of the following type
\begin{equation}\label{zN}
z_N(x) = \sum_{n=0}^{N-1} \zeta_{n,N} \cR_n(x) \approx y(x)
\end{equation}
where the coefficients $\zeta_{n,N}$ and the numerical eigenvalue $\lambda^{(N)}$ are determined by imposing, see (\ref{slp}), 
\begin{equation}\label{wf}
\sum_{n=0}^{N-1} \left\langle \cR_{m},-\cR_n'' + (q-\lambda^{(N)}) \cR_n\right\rangle \zeta_{n,N}=0, \quad \mbox{for each } m=0,\ldots,N-1.
\end{equation}
Here $\langle \cdot,\cdot \rangle$ is the standard inner product in $L_2([-1,1]),$ i.e.
$$ \langle u, v\rangle = \int_{-1}^1 u(x)v(x) dx, \qquad u,v \in L_2([-1,1]),$$
which is naturally suggested by the Liouville normal form of the SLP we are studying.
We can write (\ref{wf})  as the following generalized eigenvalue problem
\begin{equation}\label{geneig}
\left(A_N + Q_N\right)\bfzeta_N = \lambda^{(N)} B_N \bfzeta_N  
\end{equation}
where $\bfzeta_N = \left(\zeta_{0N}, \ldots,\zeta_{N-1,N}\right)^T,$ 
\begin{equation}\label{ABQ}
A_N = \left(a_{mn}\right), \quad B_N = \left(b_{mn}\right), \quad 
Q_N = \left(q_{mn}\right), \quad m,n=0,\ldots,N-1,
\end{equation}
with
\begin{equation} \label{abmn}
a_{mn} =  -\langle \cR_m, \cR_n''\rangle, \quad b_{mn} = \langle \cR_m, \cR_n\rangle, \quad 
q_{mn} = \langle \cR_m, q\,\cR_n\rangle . 
\end{equation}
The matrices $B_N$ and $Q_N$ are clearly symmetric and the same property holds for $A_N$ thanks to
the well-known Green's identity, \cite{Mag}. \\

The basis function $\cR_n$ is chosen as follows \cite{Shen94}
\begin{equation}\label{Rnp2}
 \cR_{n}(x) = \xi_{n} \cP_n(x) +  \eta_{n} \cP_{n+1}(x) +  \theta_{n} \cP_{n+2}(x) 
\end{equation}
where $\cP_j$ is the Legendre polynomial of degree $j$ 
and the three coefficients $\xi_{n}, \eta_n$ and $\theta_n$ are such that $\cR_n$ verifies the boundary conditions (BCs),
see (\ref{SN})-\ref{spanSN}).
The complete discussion of the computation of such coefficients can be found in \cite{Mag} where we used the fact that 
\begin{equation}\label{valP}
\cP_j(1) = (-1)^j \cP_j(-1) = 1, \qquad \cP'_j(1) = (-1)^{j-1}\cP'_j(-1) = j(j+1)/2
\end{equation}
and we decided to use basis functions that verify $\|\cR_n\|_\infty \leq 3$ for each $n\in \Nn_0.$  This is obtained by imposing 
$\|\left(\xi_n,\eta_n,\theta_n\right)^T\|_\infty=1$ and $\xi_n \ge 0.$ Here we simply list in
Table~\ref{tabR} the three coefficients of the linear combination in (\ref{Rnp2}) for the four problems 
subject to natural BCs and for two general ones of Robin type.\\
\begin{table}
 \caption{Coefficients $\xi_{n}, \eta_{n}$ and $\theta_{n}$ for some BCs.}\label{tabR}
\begin{center}
\begin{tabular}{|c|c||c|c|c|}
\hline
&&&&\vspace{-2mm}\\
BCs  & &$\xi_{n}$ & $\eta_{n}$ & $\theta_{n}$\\
&&&&\vspace{-2mm}\\
\hline
\hline
&&&&\vspace{-2mm}\\
$y(\pm 1)=0$ & $n\ge 0$& $1$ &  $0$ &  $-1$\\
&&&&\vspace{-2mm}\\
\hline
&&&&\vspace{-2mm}\\
$y'(\pm 1)=0$ & $n\ge 0$& $1$ & $0$ & $-\,\frac{n(n+1)}{(n+2)(n+3)}$\\
&&&&\vspace{-2mm}\\
\hline
&&&&\vspace{-2mm}\\
$\begin{array}{l} y(-1)=0 \\ \,\,\,y'(1)=0\end{array}$ & $n\ge 0$
&  $1$ & $\frac{(2n+3)}{(n+2)^2}$ &  $-\,\left(\frac{n+1}{n+2}\right)^2$\\
&&&&\vspace{-2mm}\\
\hline
&&&&\vspace{-2mm}\\
$\begin{array}{l} y'(-1)=0 \\ \quad y(1)=0\end{array}$ & $n\ge 0$
&  $1$ & $-\,\frac{(2n+3)}{(n+2)^2}$ &  $-\,\left(\frac{n+1}{n+2}\right)^2$\\
&&&&\vspace{-2mm}\\
\hline
&&&&\vspace{-2mm}\\
$\begin{array}{l}
 y'(-1) \,= \,y(-1)\\
 \,\,\,\,y'(1) \,=\, 0
\end{array}$
& $n\ge 0$
& $1$ & $\frac{2(2n+3)}{(n+2)^2(2+(n+1)(n+3))}$ & $-\,\frac{(n+1)^2(2+n(n+2))}{(n+2)^2(2+(n+1)(n+3))}$ \\
\hline
&&&&\vspace{-2mm}\\
& $n=0$ & $\frac{2}{3}$ & $1$ & $\frac{1}{3}$\vspace{-1mm}\\
$y'(\pm 1) = y(\pm 1)$ &&&&\vspace{-1mm}\\
& $n\ge 1$ & $1$ & $\frac{4(2n+3)}{(n+1)(n+2)^2(n+3)-4}$ & $-\,\frac{n(n+1)^2(n+2)+4}{(n+1)(n+2)^2(n+3)-4}$\\
\hline
\end{tabular}
\end{center}
\end{table}

For later reference, it is important to underline  the fact that, as soon as $n$ is sufficiently large, 
we always got
\begin{eqnarray}
\xi_n &=& 1,\label{xias}\\
\theta_n &=& -1 + O(n^{-1}),\label{thas}\\
\label{etas}
\eta_n &=& \left\{\begin{array}{ll}
                0, & \mbox{ if $\alpha_L\beta_R +\alpha_R\beta_L = 0,$}\\
                O(n^{-1}) & \mbox{ if $\alpha_L\beta_R+\alpha_R\beta_L \neq 0$ and $\beta_L\beta_R =0,$}\\
                O(n^{-3}) & \mbox{ if $\alpha_L\beta_R+\alpha_R\beta_L \neq 0$ and $\beta_L\beta_R \neq 0.$}
                \end{array}\right.          
\end{eqnarray}
More precisely, if the BCs are symmetric,  i.e. if $\alpha_L\beta_R + \alpha_R\beta_L =0,$ 
then we always set $\eta_n=0$ so that $\cR_n$ is an even or an odd function if $n$ is even or odd, 
respectively.

\subsection{The matrices $A_N$ and $B_N$}
In this section, we recall the results obtained in \cite{Mag} about the entries of $A_N$ and $B_N$ in (\ref{ABQ})-(\ref{abmn}).\\
Concerning the first matrix, one immediately gets that $a_{mn}=0$ for each $m>n$
since $\cR_m$ is orthogonal to any polynomial in $\Pi_{m-1},$ see (\ref{Rnp2}). Consequently,  
$A_N=A_N^T$ is diagonal with diagonal entries
\begin{eqnarray}\nonumber
a_{nn} &=&  -\xi_n\theta_n \langle \cP_n,\cP_{n+2}''\rangle \nonumber\\
    &=& -\xi_n \theta_n\left[\cP_n(x)\cP_{n+2}'(x) - \cP_n'(x)\cP_{n+2}(x) \right]_{-1}^1  \nonumber\\
    &=& - 2 (2n+3) \xi_n\theta_n,\label{amm1}
\end{eqnarray}
see (\ref{valP}). We remark that, independently of the BCs, $a_{nn}$ satisfies 
\begin{equation}\label{annas}
a_{nn} = 4\left(n+\frac{3}{2}\right) \left(1 + O\left(n^{-1}\right)\right), \quad n \gg 1.  
\end{equation}

Regarding $B_N,$ it is not too difficult to verify that it is pentadiagonal. In more detail, if we let 
$$ \hat{b}_{n}  = \langle \cP_n,\cP_n\rangle = 2/(2n+1),$$
\begin{equation} \label{R}
 \hat{B}_N = \left(\begin{array}{ccc} \hat{b}_{0}\\ &\ddots \\ &&\hat{b}_{N+1}\end{array}\right), \quad
  R_N = \left(\begin{array}{ccccc}
       \xi_0\\
       \eta_0 & \ddots\\
       \theta_0 & \ddots & \ddots\\
       & \ddots & \ddots & \xi_{N-1}\\
       && \ddots & \eta_{N-1}\\
       &&& \theta_{N-1}
\end{array}\right),
\end{equation}
then we get
\begin{equation}\label{Bfact}
B_N = R_N^T \, \hat{B}_N R_N.
\end{equation}

\subsection{The matrix $Q_N.$}\label{secQ}
From (\ref{ABQ})-(\ref{abmn}), one obtains that $Q_N$
admits a factorization similar to the one just given for $B_N.$ Specifically
\begin{eqnarray} \label{Qfact}
Q_N &=& R_N^T \, \hat{Q}_N R_N, \qquad 
\hat{Q}_N = \left(\hat{q}_{mn}\right) \label{hQ}
 \end{eqnarray}
where $R_N$ is defined in (\ref{R}),
\begin{equation}
\hat{q}_{mn} =  \langle \cP_m,q\cP_n\rangle  
= \sum_{i=1}^S \langle \cP_m,g_i\cP_n\rangle_{(\be_i,\ga_i)}
\equiv  \sum_{i=1}^S \hat{q}_{mn}^{(i)},\label{hqmn1}
\end{equation}
being
\begin{equation}\label{pscalbg}
\langle u, v\rangle_{(\be,\ga)} = \int_{-1}^1 \frac{u(x)\,v(x)}{(1-x)^\be\,(1+x)^\ga} dx, \qquad \be,\ga<1.
\end{equation}

As done in \cite{GM,Mag}, it is possible to prove the following result 
by using the well-known recurrence relation of the Legendre polynomials.

\begin{proposition}\label{propGFinf} Let $q\in L_1([-1,1])$ and 
\begin{eqnarray*}
\hat{\bfq}_n &\equiv& \left(\hat{q}_{0n}, \hat{q}_{1n}, \ldots \right )^T  \in \ell_\infty, \quad n \ge 0.
\end{eqnarray*}
If we define the linear tridiagonal operator $ \bfz \in \ell_\infty \mapsto \mathcal{H}\,\bfz \in \ell_\infty$ where
\begin{eqnarray}\nonumber
\mathcal{H} &=&
\left( \begin{array}{cccccc}
            0 & h_{01}\\
            h_{10} & 0 & h_{12}\\
            &h_{21} & 0 & h_{23}\\
            && \ddots & \ddots & \ddots\\
            \end{array}\right), \quad \begin{array}{l} h_{m,m-1}=m/(2m+1), \\ \\ h_{m,m+1}=(m+1)/(2m+1),\end{array}    \label{H}
\end{eqnarray}
and we let $\hat{\bfq}_{-1}$ be the zero sequence, then we get 
\begin{equation} \label{recvn}
\hat{\bfq}_{n+1} = \frac{2n+1}{n+1}\, \mathcal{H}\hat{\bfq}_n -\frac{n}{n+1} \,\hat{\bfq}_{n-1},\qquad  n\ge 0. \qquad \mbox{\qed}\\
\end{equation}
\end{proposition}

Now, the structure of $\mathcal{H}$ and (\ref{recvn}) permit to 
determine the entire matrix $\hat{Q}_N$ once $\hat{q}_{m0}$ 
have been computed for each $m=0,1,\ldots,2N+2$ (see \cite{GM,Mag} for the details). We shall proceed  
by discussing how we determine these values for problems with $S=1$ since the generalization is simple, see (\ref{hqmn1}).
In this regard, we observe that (\ref{pscalbg}), \cite[16.4 formula (2)]{EMOT} and arguments similar 
to the ones used in the proof of \cite[Proposition 2]{GM}
allow to get that
\begin{equation}\label{hq0}
\hat{\bfq}_0 = g_1(\mathcal{H}) \left(\begin{array}{c}\hat{q}_0^{(1)}\\ \hat{q}_1^{(1)}\\\vdots \end{array}\right), \quad
\hat{q}_m^{(1)} = \langle \cP_m,\cP_0\rangle_{(\be_1,\ga_1)}=\langle \cP_m,1\rangle_{(\be_1,\ga_1)}. 
\end{equation}
Let us assume for the moment that we have computed the values of $\hat{q}_m^{(1)}$ for each $m=0,\ldots,L$
with $L$ sufficiently large. Then, recalling that by assumption $g_1$ is analytical inside and over a Bernstein ellipse 
containing $[-1,1]$, we proceed in this way. We get a polynomial approximation of $g_1$ by transforming it in a 
{\tt Chebfun} function \cite{Cf}, which is accurate up to machine precision, and then we apply 
the previous formula to compute the first $2N+2$ entries of $\hat{\bfq}_0.$\\
Concerning the computation of $\hat{q}_m^{(1)}$ 
we have to distinguish the following cases:
\begin{enumerate}
 \item $\be_1=\ga_1=0:$ it is evident that $\hat{q}_{0}^{(1)}=2$ and $\hat{q}_{m}^{(1)}=0$ for each $m>0;$
 \item $\be_1=0, \ga_1\neq 0:$ as discussed in \cite{Mag} it results
 \begin{equation}\label{be1u0}
  \hat{q}_m^{(1)} =\frac{(-1)^m\,2^{1-\ga_1}\,(\ga_1)_m}{(1-\ga_1)_{m+1}}  
 \end{equation} 
 where $(t)_\ell$ is the Pochhammer symbol;
 \item $\be_1\neq 0, \ga_1= 0:$ with similar computations one gets
 \begin{equation}\label{ga1u0}
 \hat{q}_m^{(1)} =\frac{2^{1-\be_1}\,(\be_1)_m}{(1-\be_1)_{m+1}};
 \end{equation}
 \item $\be_1\ga_1\neq 0:$ by using \cite[16.2, formula(6)]{EMOT} we get
 \begin{equation}\label{gabenoz}
  \hat{q}_m^{(1)} = \alpha
                    \,_3F_2\left(\begin{array}{c} -m, \,\, 1+m,\,\,1-\be_1\\ 1,\,\, 2-\be_1-\ga_1\end{array};1\right)
 \end{equation}
 where
 \begin{equation}\label{ombega}
  \alpha = \frac{2^{1-\be_1-\ga_1} \Gamma(1-\be_1)\Gamma(1-\ga_1)}{\Gamma(2-\be_1-\ga_1)} = \hat{q}_0^{(1)}.
 \end{equation}

\end{enumerate}
It is clear that  the formulas in (\ref{be1u0}) or (\ref{ga1u0}) allow to compute $\hat{q}_m^{(1)}$ easily.
For example, if $\be_1\neq 0$ and  $\ga_1= 0$ then from (\ref{ga1u0}) one gets 
$$ \hat{q}_0^{(1)} = \frac{2^{1-\be_1}}{1-\be_1}, \qquad \quad {q}_{m+1}^{(1)}= \frac{m+\be_1}{m+2-\be_1}{{q}_{m}^{(1)}}, \quad m\ge 0,$$
so that it is possible to proceed recursively.
On the other hand, if $\be_1\ga_1\neq 0$ then the computation of the Gauss hypergeometric function 
at right hand-side of (\ref{gabenoz}) can be costly and ill-conditioned. 
We thus preferred to find alternative expressions.  
In particular, if $\be_1=\ga_1\neq 0$ then the application of the following Whipple sum
\begin{eqnarray*}
\lefteqn{\,_3F_2\left(\begin{array}{c} a, \,\, 1-a,\,\,c\\ e,\,\, 1+2c-e\end{array};1\right) =}\\
&&\frac{2^{1-2c}\pi \Gamma(e) \Gamma(1+2c-e)}{\Gamma((a+e)/2) \Gamma((a+1+2c-e)/2) \Gamma((e+1-a)/2)\Gamma((2+2c-a-e)/2)} 
\end{eqnarray*}
with $a=-m,\, c=1-\ga_1$ and $e=1$ gives
$$ \hat{q}_m^{(1)} = \frac{\pi \Gamma^2(1-\ga_1)}{\Gamma\left((3-2\ga_1 +m)/2\right)\Gamma\left((2-2\ga_1-m)/2\right) \Gamma\left((2+m)/2\right)\Gamma\left((1-m)/2\right)}.$$
Therefore, if $m$ is odd then $\hat{q}_m^{(1)}=0$ (this was indeed already evident from its definition in (\ref{hq0}) with $\be_1=\ga_1$). On the other hand, 
if $m$ is even then
\begin{eqnarray*}
\hat{q}_0^{(1)}     &=& \frac{\sqrt{\pi}\,\, \Gamma(1-\ga_1)}{\Gamma(3/2-\ga_1)} = \frac{2^{1-2\ga_1}\,\Gamma^2(1-\ga_1)}{\Gamma(2-2\ga_1)}, \\
\hat{q}_{m+2}^{(1)} &=& \frac{(m+2\ga_1)(m+1)}{(m+3-2\ga_1)(m+2)}\,\hat{q}_m^{(1)}, \qquad m=0,2,4,\ldots\,.
\end{eqnarray*}

\begin{remark}\label{simm}
If $q(x) = q(-x)$ and $N$ is even then the matrix $\hat{Q}_N$ in (\ref{hQ}) is permutation similar to a $2 \times 2$ 
block diagonal matrix with diagonal blocks of size $N/2+1.$
In addition, if the coefficients of the BCs in (\ref{bca})-(\ref{bcb}) verify $\alpha_L\beta_R+\alpha_R\beta_R=0$
then the SLP is called symmetric and, as we recalled after (\ref{etas}), $\eta_n=0$ for each $n.$
Therefore, the matrices $A_N,$ $Q_N$ and $B_N$ are 
permutation similar to $2 \times 2$ block diagonal matrices too, see the paragraph before (\ref{amm1}) and
(\ref{R})--(\ref{Qfact}). This implies that we can split the generalized eigenvalue problem (\ref{geneig}) 
into two ones of halved size.
\end{remark}
It remains to discuss how it is possible to avoid the evaluation  of the Gauss hypergeometric function for
the computation of $\hat{q}_m^{(1)}$ if  $\be_1\ga_1\neq 0$ and $\be_1\neq \ga_1.$  In this case, even though alternative strategies are possible, 
we decided to write $\cP_m\equiv \cP^{(0,0)}_m$ as a linear combination of $\left\{\cP_\ell^{(0,-\ga_1)}\right\}_{\ell=0}^m$ where $\cP_\ell^{(0,-\ga_1)}$ is 
the Jacobi polynomial of degree $\ell$ with weighting function $\omega(x) =(1+x)^{-\ga_1}.$ In other words, first of all we write
$$\cP_m(x) \equiv \cP_m^{(0,0)}(x) = \sum_{\ell=0}^{m} \chi_{m,\ell} \cP_\ell^{(0,-\ga_1)}(x).$$
Then we determine $\hat{q}_m^{(1)}$ as follows
$$\hat{q}_m^{(1)} = \sum_{\ell=0}^m \chi_{m,\ell} \,\langle \cP_\ell^{(0,-\ga_1)},1\rangle_{(\be_1,\ga_1)}.$$ 
Now, by using formulas in \cite[16.4]{EMOT} we obtain
\begin{eqnarray*}
\chi_{m,\ell} &=& \frac{\langle\cP_m^{(0,0)},\cP_\ell^{(0,-\ga_1)}\rangle_{(0,\ga_1)}}
                       {\langle\cP_\ell^{(0,-\ga_1)},\cP_\ell^{(0,-\ga_1)}\rangle_{(0,\ga_1)}}\\
              &=&     
                   \frac{\Gamma(1-\ga_1) (2\ell+1-\ga_1)  \Gamma(m+\ell+1) \Gamma(\ell+1-\ga_1)}
                        {\Gamma(\ell+1) \,\Gamma(m+\ell+2-\ga_1) \Gamma(m-\ell+1) \Gamma(\ell-m+1-\ga_1)},\\
\end{eqnarray*}
$$
\langle \cP_\ell^{(0,-\ga_1)},1\rangle_{(\be_1,\ga_1)} =
                  \frac{2^{1-\be_1-\ga_1} \Gamma(1-\be_1) \Gamma(\ell+1-\ga_1) (\be_1)_\ell}
                       {\Gamma(2-\be_1-\ga_1)\,\Gamma(\ell+1)\, (2-\be_1-\ga_1)_\ell}.
$$                       
Hence
$$ \hat{q}_m^{(1)} = \alpha \sum_{\ell=0}^m t_{m-\ell} u_{m+\ell} \nu_\ell, $$
where $\alpha$ is defined in (\ref{ombega}),
\begin{eqnarray*}
 t_r &=& \frac{1}{\Gamma(1+r)\Gamma(1-\ga_1-r)},\qquad 
 u_r = \frac{\Gamma(r+1)}{\Gamma(r+2-\ga_1)}, \label{trur}\\
 \nu_\ell &=& (2\ell+1-\ga_1) \left(\frac{\Gamma(\ell+1-\ga_1)}{\Gamma(\ell+1)}\right)^2\,\frac{(\be_1)_\ell}{(2-\be_1-\ga_1)_\ell}.
 \label{nul}
\end{eqnarray*}
It is not too difficult to verify that 
\begin{equation}\label{qm1vet}
\left(
\begin{array}{c}
\hat{q}_0^{(1)}\\
\hat{q}_1^{(1)}\\
\hat{q}_2^{(1)}\\
\hat{q}_3^{(1)}\\
\vdots
\end{array}
\right) = \alpha
\left( T \circ U \right) 
\left(
\begin{array}{c}
\nu_0\\
\nu_1\\
\nu_2\\
\nu_3\\
\vdots
\end{array}
\right),
\end{equation}
being ``$\circ$'' the Hadamard product and 
$$ T = \left(\begin{array}{ccccc}
             t_0 \\
             t_{1} & t_0 \\
             t_{2} & t_{1} & t_0\\
             t_3 & t_{2} & t_{1} & t_0\\
             \vdots & \ddots & \ddots & \ddots & \ddots               
             \end{array}\right),
\qquad 
 U = \left(\begin{array}{ccccc}
           u_0 & u_1 & u_2 & u_3 &\ldots\\
           u_1 & u_2 & u_3 \\
           u_2 & u_3\\
           u_3\\
           \vdots
           \end{array}\right),
$$
i.e. $T$ and $U$ are a lower triangular Toeplitz and an Hankel matrix, respectively. 
Clearly, a suitable truncation of the vectors and matrices in (\ref{qm1vet}) is operated depending on the 
number of values of $\hat{q}_m^{(1)}$ that we actually need. In addition, we compute $\nu_\ell,$ $t_r$ and $u_r$ recursively.
Finally, it is worth mentioning that an algorithm similar to the one described in \cite{TWO}
can be used for the matrix-vector product in (\ref{qm1vet}).

\section{Error analysis and computation of corrected numerical eigenvalues.} \label{sec3}
We now  study the behavior of the error in the resulting numerical eigenvalues 
as $N$ increases and for a fixed index. 
In particular, we consider weakly regular problems with a potential of the type specified in (\ref{pot1}) which is 
unbounded  at least at one endpoint. 
The analysis that we are going to present will be also used to derive a very effective and 
efficient procedure for an a posteriori correction of the 
numerical eigenvalues.\\

Let $\lambda^{(N)}$ be the approximation of the exact eigenvalue $\lambda$ as $N$ increases and let  
$y$ be the corresponding exact eigenfunction having the following expansion
\begin{equation}
 y(x) = \sum_{n=0}^{+\infty} c_n \cR_n(x). \label{espandiy} 
\end{equation}

The following first result can be proved by using arguments similar to the ones considered in \cite{Mag}.

\begin{theorem}
If $N$ is sufficiently larger than the index of the eigenvalue then, see (\ref{zN}) and (\ref{amm1}),
\begin{eqnarray}
 \nonumber c_n &\approx& -\,\frac{\langle \cR_n, q y\rangle}{a_{nn}} = -\,\sum_{i=1}^S \frac{\langle \cR_n, g_i y\rangle_{(\be_i,\ga_i)}}{a_{nn}},\\
\lambda-\lambda^{(N)} &\approx& - \frac{1}{\langle z_N,y\rangle} \sum_{n=N}^{+\infty} c_n \,\langle \cR_n, qz_N\rangle
\approx -\frac{1}{\langle z_N,y\rangle}\sum_{i,j=1}^S \Delta_{ij} \label{error3} 
\end{eqnarray}
where 
\begin{equation}
\label{Deij}
\Delta_{ij} = \sum_{n=N}^{+\infty} \frac{1}{a_{nn}}\langle \cR_n, g_i y\rangle_{(\be_i,\ga_i)}\,\langle \cR_n, g_j z_N\rangle_{(\be_j,\ga_j)}.
                                         \quad \mbox{\qed}
\end{equation}
\end{theorem}

The asymptotic estimate that we are going to prove in the next theorem is fundamental for proceeding.
\begin{theorem}\label{propsi} Let $\psi \in C^{\infty}(-1,1)\bigcap C^1[-1,1].$ 
If $\psi(-1)\,\psi(1)\neq 0,$ $\be,\ga <1$ and if $n$ is sufficiently large then
\begin{equation}\label{sidi}
\langle \cR_n,\psi\rangle_{(\be,\ga)} \approx 
\frac{(-1)^n 2\,\omega(-1,\ga,\be)\, \psi(-1)}{(n+3/2)^{p(\ga)/2}} + \frac{2\,\omega(+1,\be,\ga)\, \psi(1)}{(n+3/2)^{p(\be)/2}},
\end{equation}
where
\begin{equation}
p(\delta) = 6-4\delta, \qquad
\omega(\pm 1,\delta_0,\delta_1) = 2(2-\delta_0-\kappa_{\pm 1})\,\hat{\omega}(\delta_0,\delta_1), \label{pdomxdd}
\end{equation}
being
\begin{equation} \label{kappa}
\kappa_{\pm 1} = \left\{\begin{array}{ll} 1, & \mbox{  if $\cR_n(\pm 1)=0$}\\
                                          0, & \mbox{  otherwise} 
                      \end{array}\right., \quad
\hat{\omega}(\delta_0,\delta_1) = \frac{2^{1-\delta_0-\delta_1}\Gamma(1-\delta_0)}{\Gamma(\delta_0)}.
\end{equation}              
\end{theorem}
\underline{Proof}: Recalling the definition of $\cR_n$ in (\ref{Rnp2}), let us consider
first of all $\langle \cP_n,\psi\rangle_{(\be,\ga)}.$
In this regard, if we use the results proved in \cite{Sid09}, (\ref{hq0})--(\ref{ga1u0}) and we assume that
$n$ is sufficiently large then we get
\begin{eqnarray*}
\langle \cP_n,\psi\rangle_{(\be,\ga)} &\approx&
\frac{\psi(-1)}{2^{\be}} \langle \cP_n,1\rangle_{(0,\ga)}  + 
\frac{\psi(1)}{2^{\ga}} \langle \cP_n,1\rangle_{(\be,0)}\\
&=& \frac{\psi(-1)\,(-1)^n\, (\ga)_n}{2^{\be+\ga-1} (1-\ga)_{n+1}} +
 \frac{\psi(1) \,(\be)_n}{2^{\be+\ga-1} (1-\be)_{n+1}}\\
&=& \frac{\psi(-1)\, (-1)^n\,\hat{\omega}(\ga,\be)\,\Gamma(n+\ga)}{\Gamma(n+2-\ga)}+
 \frac{\psi(1)\,\hat{\omega}(\be,\ga) \,\Gamma(n+\be)}{\Gamma(n+2-\be)}.
\end{eqnarray*}
Now, it is known that the ratio of two gamma functions satisfies
$$\frac{\GG(z+a)}{\GG(z+b)} = z^{a-b} \left(1 + \frac{(a-b)(a+b-1)}{2z} + O (z^{-2})\right),\quad z\gg 0.$$
Therefore, if we use it with $z = n+1/2$ then we obtain
$$
\langle \cP_n,\psi\rangle_{(\be,\ga)} \approx
  \left(\frac{\psi(-1)\,(-1)^n\,\hat{\omega}(\ga,\be)}{(n+1/2)^{\hat{p}(\ga)}}+
\frac{\psi(1)\,\hat{\omega}(\be,\ga)}{(n+1/2)^{\hat{p}(\be)}}\right)\left(1+O(n^{-2})\right),  
$$
with
\begin{equation}\label{pLR}
\hat{p}(\delta) = 2-2\delta= \frac{p(\delta)}2-1, \qquad \delta=\ga,\be.
\end{equation}
This implies that to determine an estimate for $\langle \cR_n,\psi\rangle_{(\be,\ga)}$, we have to study
these terms 
\begin{eqnarray*}
\lefteqn{(-1)^n \left( \frac{\xi_n}{(n+1/2)^{\hat{p}(\ga)}} -  \frac{\eta_n}{(n+3/2)^{\hat{p}(\ga)}} + \frac{\theta_n}{(n+5/2)^{\hat{p}(\ga)}}\right)}\\
 && \qquad \qquad \approx (-1)^n\left(n+\frac{3}{2}\right)^{-\hat{p}(\ga)} \left(\xi_n-\eta_n+\theta_n  + \frac{\hat{p}(\ga)(\xi_n-\theta_n)}{n+3/2}\right),\\
\lefteqn{\frac{\xi_n}{(n+1/2)^{\hat{p}(\be)}} +  \frac{\eta_n}{(n+3/2)^{\hat{p}(\be)}} + \frac{\theta_n}{(n+5/2)^{\hat{p}(\be)}}}\\
&& \qquad \qquad \approx \left(n+\frac{3}{2}\right)^{-\hat{p}(\be)} \left(\xi_n+\eta_n+\theta_n  + \frac{\hat{p}(\be)(\xi_n-\theta_n)}{n+3/2}\right).
\end{eqnarray*}

We recall that if $n$ is sufficiently large then $\xi_n = 1,$ see (\ref{xias}). In addition,
by using the formulas in \cite{Mag}, see also (\ref{thas})-(\ref{etas}), it is possible to verify with some computations that 
\begin{enumerate}
 \item if $\cR_n(-1) =(-1)^n (\xi_n-\eta_n+\theta_n) \neq 0$  then
  $$ \xi_n-\eta_n+\theta_n = \frac{4}{n+3/2} \left(1 + O\left(\frac{1}{n}\right)\right);$$
 \item if $\cR_n(1) =\xi_n+\eta_n+\theta_n\neq 0$ then
  $$ \xi_n+\eta_n+\theta_n = \frac{4}{n+3/2} \left(1 + O\left(\frac{1}{n}\right)\right);$$
 \item $\xi_n -\theta_n = 2 \left(1 + O\left(n^{-1}\right)\right).$
  \end{enumerate}
Therefore, see (\ref{kappa}) and (\ref{pLR}), 
\begin{eqnarray*}
\xi_n-\eta_n+\theta_n  + \frac{\hat{p}(\ga)(\xi_n-\theta_n)}{n+3/2} &\approx& \frac{4(2-\ga-\kappa_{-1})}{n+3/2}\\
\xi_n+\eta_n+\theta_n  + \frac{\hat{p}(\be)(\xi_n-\theta_n)}{n+3/2} &\approx& \frac{4(2-\be-\kappa_{+1})}{n+3/2}.
\end{eqnarray*}
The statement follows by collecting all these partial results.\qed\\

It must be underlined that $\hat{\omega}(\delta_0,\delta_1)=0$ if $-\delta_0 \in \Nn_0,$ see (\ref{kappa}).
This implies that one or both the terms at the right hand-side of (\ref{sidi}) can be zero.
For our purposes, this does not constitute a problem since in the convergence analysis that we are going to
prove such terms are surely negligible with respect to the others.\\

We need the following notation to proceed: for each $i=1,\ldots,S,$ let
\begin{equation}\label{rili}
g_i(x) = (1-x)^{r_{i}}(1+x)^{\ell_{i}}\hat{g}_i(x),  \quad  \hat{g}_i(\pm 1) \neq 0,
\end{equation}
i.e. let $r_i$ and $\ell_i$ be the multiplicities of $x=1$ and $x=-1,$ as zeros of $g_i,$ respectively.
We are now ready for proving the following theorem.
\begin{theorem}[Convergence]\label{conv}
Let assume the potential in (\ref{pot1}) is unbounded at least at one endpoint and let consider the 
following subsets of $\left\{1,2,\ldots,S\right\}$
\begin{eqnarray*}
{\cal{I}}_L &\subseteq& \left\{1,2,\ldots,S\right\} \qquad \mbox{with} \qquad i \in  {\cal{I}}_L \quad \longleftrightarrow \quad -\gamma_i \notin \Nn_0,\\
{\cal{I}}_R &\subseteq& \left\{1,2,\ldots,S\right\} \qquad \mbox{with} \qquad i \in  {\cal{I}}_R \quad \longleftrightarrow \quad -\beta_i \notin \Nn_0.
\end{eqnarray*}
If $N$ is sufficiently larger than the index of the eigenvalue then
\begin{equation}\label{ord}
\lambda-\lambda^{(N)} = O\left((N+1)^{-p}\right), \qquad p =\inf\left(p_L,p_R\right),
\end{equation}
where, see (\ref{pdomxdd})-(\ref{kappa}) and (\ref{rili}),
\begin{eqnarray*}
p_L &=& \inf_{i\in {\cal I}_L} p(\hat{\ga}_i)=\inf_{i\in {\cal I}_L} 6-4\hat{\ga}_i, \qquad \hat{\ga}_i=\ga_i-\ell_i-\kappa_-,\\
p_R &=& \inf_{i\in {\cal I}_R} p(\hat{\be}_i)=\inf_{i\in {\cal I}_R} 6-4\hat{\be}_i, \qquad \hat{\be}_i=\be_i-\ell_i-\kappa_+,
\end{eqnarray*}
being $\inf\emptyset = +\infty$ by convention.
\end{theorem}
\underline{Proof}: 
From the definition of $\kappa_{\pm},$  (\ref{zN}) and (\ref{espandiy}) we get  
\begin{eqnarray}
y(x) &=& (1-x)^{\kappa_+}(1+x)^{\kappa_-} \hat{y}(x), \qquad \quad \hat{y}(\pm 1) \neq 0, \label{hy}\\
z_N(x) &=& (1-x)^{\kappa_+}(1+x)^{\kappa_-} \hat{z}_N(x), \qquad \hat{z_N}(\pm 1) \neq 0,\label{hz}
\end{eqnarray}
so that recalling (\ref{error3})-(\ref{Deij})  we must determine an estimate for 
$$ \frac{1}{a_{nn}}\left(\left\langle \cR_n, \hat{g}_i \hat{y}\right\rangle_{(\hat{\be}_i,\hat{\ga}_i)}
\left\langle \cR_n, \hat{g}_j \hat{z}_N\right\rangle_{(\hat{\be}_j,\hat{\ga}_j)}\right) \equiv \left(\star\right).$$
To this end, we apply (\ref{sidi}) with $(\be,\ga)=(\hat{\be}_i,\hat{\ga}_i)$ and $\psi=\hat{g}_i \hat{y},$  
or with $(\be,\ga)=(\hat{\be}_j,\hat{\ga}_j)$ and $\psi=\hat{g}_j \hat{z}_N.$ In this way, recalling also (\ref{annas}), we obtain
\begin{eqnarray*}
\left(\star\right)&\approx& \quad \,\left(\frac{(-1)^n \sigma_{iL}\,\hat{y}(-1)}{(n+3/2)^{1+p(\hat{\ga}_i)/2}} + 
                    \frac{\sigma_{iR}\,\hat{y}(1)}{(n+3/2)^{1+p(\hat{\be}_i)/2}}\right)\\
&& \qquad \qquad \,\,\times 
\left(\frac{(-1)^n \sigma_{jL}\,\hat{z}_N(-1)}{(n+3/2)^{p(\hat{\ga}_j)/2}} + 
                    \frac{\sigma_{jR}\,\hat{z}_N(1)}{(n+3/2)^{p(\hat{\be}_j)/2}}\right),
\end{eqnarray*}
where 
$$
\sigma_{iL} = \hat{g}_i(-1) \omega(-1,\hat{\ga}_i,\hat{\be}_i),\qquad
\sigma_{iR} = \hat{g}_i(+1) \omega(+1,\hat{\be}_i,\hat{\ga}_i).
$$
We now  use the following integral estimates with suitable $\bar{p}>0:$ 
\begin{eqnarray*}
\sum_{n=N}^{+\infty} \frac{1}{(n+3/2)^{1+\bar{p}}} &\approx& \int_{N}^{+\infty} (n+1)^{-1-\bar{p}} dn = \frac{1}{\bar{p} (N+1)^{\bar{p}}},\\
\sum_{n=N}^{+\infty} \frac{(-1)^n}{(n+3/2)^{1+\bar{p}}}  
&{\begin{array}{c} \ell =\mbox{rem}(N,2)\\ \uparrow\\ \approx\\~\\\end{array}}& (-1)^N \sum_{m=(N-\ell)/2}^{+\infty} \frac{1+\bar{p}}{(2m+\ell+3/2)^{2+\bar{p}}}\\
&\approx& (-1)^N\int_{(N-\ell)/2}^{+\infty} \frac{(1+\bar{p}) \,\, dm}{(2m+\ell+1)^{2+\bar{p}}} \\
&=& \frac{(-1)^N}{2 (N+1)^{1+\bar{p}}}.
\end{eqnarray*}
In particular, we apply the first one with $\bar{p} = (p(\hat{\ga}_i)+p(\hat{\ga}_j))/2,$ or
$\bar{p} = (p(\hat{\be}_i)+p(\hat{\be}_j))/2$ and the second estimate with $\bar{p} = (p(\hat{\ga}_i)+p(\hat{\be}_j))/2$ or 
with $\bar{p} = (p(\hat{\be}_i)+p(\hat{\ga}_j))/2.$ This leads to, see (\ref{Deij}), 
\begin{eqnarray}
\label{bDel} \Delta_{ij} \approx \bar{\Delta}_{ij} = && \frac{2\,\sigma_{iL}\,\sigma_{jL}\hat{y}(-1)\hat{z}_N(-1)}{(p(\hat{\ga}_i)+p(\hat{\ga}_j))(N+1)^{(p(\hat{\ga}_i)+p(\hat{\ga}_j))/2}}\\
\nonumber           &+&\frac{2\,\sigma_{iR}\,\sigma_{jR}\,\hat{y}(1)\hat{z}_N(1)}{(p(\hat{\be}_i)+p(\hat{\be}_j))(N+1)^{(p(\hat{\be}_i)+p(\hat{\be}_j))/2}}\\
\nonumber           &+&  \frac{(-1)^N\,\sigma_{iL}\,\sigma_{jR}\,\hat{y}(-1)\hat{z}_N(1)}{2(N+1)^{1+(p(\hat{\ga}_i)+p(\hat{\be}_j))/2}}\\
\nonumber           &+&  \frac{(-1)^N\,\sigma_{iR}\,\sigma_{jL}\,\hat{y}(1)\hat{z}_N(-1)}{2(N+1)^{1+(p(\hat{\be}_i)+p(\hat{\ga}_j))/2}}.
\end{eqnarray}
Therefore
\begin{equation}\label{errdb}
\lambda-\lambda^{(N)} \approx -\frac{1}{\langle z_N,y\rangle}\sum_{i,j=1}^S\bar{\Delta}_{ij} 
\end{equation}
and the statement follows by observing that the principal term of such summation behaves like $O\left((N+1)^{-p}\right)$  
where $p$ is defined in (\ref{ord}). \qed\\

As done in \cite{Mag}, we now discuss how we can use (\ref{errdb}) to improve the accuracy of the numerical 
eigenvalues. The approach is that of considering the following normalization for the numerical 
and the exact eigenfunctions, see (\ref{hy})-(\ref{hz}),   
$$
\langle z_N,z_N\rangle = \bfzeta_N^T B_N \bfzeta_N =1, \qquad \langle y,y \rangle = 1,\qquad \hat{z}_N(-1),\hat{y}(-1) >0.\\
$$
By using the orthogonality of the Legendre polynomials,  the estimates
$\langle z_N,y \rangle \approx 1$ and $\hat{y}(\pm 1) \approx \hat{z}_N(\pm 1)$  follow and
consequently the next formula for the computation of corrected numerical eigenvalues
\begin{equation}\label{mu}
\mu^{(N)} = \lambda^{(N)} - \sum_{i,j=1}^S \hat{\Delta}_{ij}
\end{equation}
where $\hat{\Delta}_{ij}$ is obtained from $\bar{\Delta}_{ij}$  via the 
substitutions $ \hat{y}(\pm 1) \rightarrow \hat{z}_N(\pm 1),$ see (\ref{hy})--(\ref{bDel}).
These are done by using $z_N(\pm 1)$ or $z_N'(\pm 1).$ 
For example, if $y(-1)=0$ and $y(1)\neq 0$ then $y(x)=(1+x)\hat{y}(x)$ and
$z_N(x)=(1+x)\hat{z}_N(x).$ Consequently
$$ \hat{y}(1) = y(1)/2 \approx z_N(1)/2=\hat{z}_N(1), \qquad 
    \hat{y}(-1) = y'(-1) \approx z_N'(-1)=\hat{z}_N'(-1).$$

Finally, we must say that in the actual implementation we do not consider 
the correction terms corresponding to values of 
$r_i$ and of $\ell_i$ in (\ref{rili}) that are greater than one since
their contributions are surely irrelevant with respect to the others.
The advantage is that for the computation of $\mu^{(N)}$
we do not need to evaluate derivatives of $g_i$ at the endpoints of order greater than one.
In this way, altogether, the cost for the application of (\ref{mu})    
is very low since it is essentially given by the evaluations of $g_i(\pm 1),$ eventually of $g_i'(\pm 1)$  and of    
$\hat{z}_N(\pm 1)$ which is simple because the values of $\cR_n(\pm 1)$ or of $\cR_n'(\pm 1)$ are known
(see (\ref{Rnp2}) and (\ref{valP})).

\section{Numerical tests} \label{sec4}
The method described was implemented in Matlab ({\tt ver.R2017a}). In particular, 
we used routines included in the open-source {\tt Chebfun} package \cite{Cf} for the 
computation of the matrix $Q_N$ and we solved the generalized eigenvalue problem (\ref{geneig})
by using the {\tt eigs}  function, with option ``{\tt SM}'' for getting the ones of smallest magnitude.\\

\begin{figure}[h]
\begin{center}
 \includegraphics[width=.9\textwidth, height = .7\textwidth]{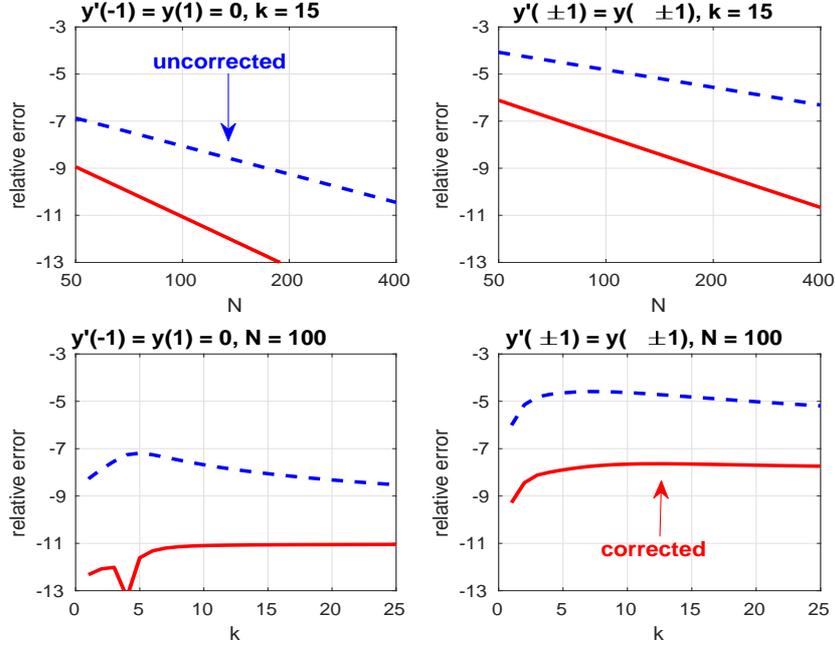}
\caption{Relative errors for problems with potential (\ref{q2}).}
 \label{fig1}
\end{center}
\end{figure}

The first results that we present confirm the statement of Theorem~\ref{conv}. 
In particular, we considered the problems with one of the following potentials
\begin{eqnarray}
q(x) &=& \frac{10 \,\mathrm{e}^{1-x}}{(1-x)^{3/4}(1+x)^{1/4}},   \label{q1}\\
q(x) &=& \frac{10\cos(4(1+x))}{(1+x)^{1/2}} + \frac{5\sin(4(1+x))}{(1-x)^{7/8}(1+x)^{3/4}}, \label{q2}
\end{eqnarray}
subject to one of the next four BCs
\begin{equation}\label{BCq12}
y'(-1)=y(1)=0, \quad y(-1)=y'(1)=0, \quad y'(\pm 1) = 0, \quad y'(\pm 1)= y(\pm 1).
\end{equation}
In addition, we used the classical formula 
$$p  \approx \log_2\left(\delta\lambda_{k,N}/\delta\lambda_{k,2N+1}\right), \qquad \delta\lambda_{k,N} \equiv |\lambda_k^{(N)}-\lambda_k^{(2N+1)}|,$$
for the numerical estimate of the order of convergence (the lower index $k$ denotes 
the index of the eigenvalue).
The results we got for the eigenvalues of 
index $k=5,10,20,$ are listed in Table~\ref{tab1} for the first potential and in Table~\ref{tab2} for
the second one. As one can see, such results are in perfect agreement with the statement of 
the theorem previously mentioned.\\
Concerning the problems with $q$ defined in (\ref{q2}) subject to the first or to the last BCs in (\ref{BCq12}),  
we applied the a posteriori correction, namely we computed also $\mu_k^{(N)}$ 
defined in (\ref{mu}). In addition, for these problems and for the subsequent ones, we evaluated the relative errors 
\begin{equation}\label{relerr}
\log_{10}\left(|\lambda_{k}^{(N)}-\bar{\lambda}_k|/|\bar{\lambda}_k|\right),
\qquad \log_{10}\left(|\mu_{k}^{(N)}-\bar{\lambda}_k|/|\bar{\lambda}_k|\right),
\end{equation}
by using as reference ``exact'' eigenvalue 
the values of $\bar{\lambda}_k\equiv \mu_k^{(Nt)}$ with $Nt \gg N\ge k.$
As discussed in the introduction, this choice was motivated by the fact that
the accuracy of the numerical approximations of $\lambda_k$ provided by the 
MATSLISE2 \cite{Led}, the SLEDGE \cite{Sledge} and the SLEIGN2 \cite{Sleign} codes
may be  not sufficient for our purposes.
The resulting relative errors (\ref{relerr}) have been reported in Figure~\ref{fig1}. In more details, 
in the subplots at the top of such figure, the relative errors in the approximation of 
the fifteenth eigenvalue are plotted versus $N$ with $N$ ranging from $50$ to $400.$
For the subplots at the bottom, instead, we fixed $N=100$ and we depict 
the errors for the index $k$ ranging from $1$ to $25.$ The legend of each graphic and of the subsequent ones 
is dashed line and solid line for the errors in the uncorrected numerical eigenvalues and in the corrected ones, respectively.
These results show that the a posteriori correction is very effective from many point of views. In fact:
\begin{itemize}
 \item from the subplots on the bottom  one deduces that for $N=100$ 
and $1\leq k \leq 25$ the gain resulting from the correction is larger than 
two significant digits;
\item the two subplots on the top show that  $\mu_{15}^{(N)}$
is always more accurate than $\lambda_{15}^{(2N)};$ 
\item the error in the corrected numerical eigenvalues decreases much faster with respect to $N$
than the error in the uncorrected ones. Concerning this point,
we used a least-square approach to evaluate numerically the
order of convergence $p_\mu$ such that 
$$|\mu_{k}^{(N)}-\bar{\lambda}_{k}| = O((N+1)^{-p_\mu}).$$
For these examples, we obtained
$p_\mu \approx 7$ for the problem subject to Neumann-Dirichlet BCs and $p_\mu \approx 5$ for the one 
subject to the unsimmetric Robin-Robin BCs.
\end{itemize}
\begin{table}[t]
\caption{Order of convergence for problems with potential (\ref{q1}).} \label{tab1}
\begin{small}
\begin{center}
\begin{tabular}{|r|c c|c c|c c|}\hline
\multicolumn{7}{|c|}{~}\vspace{-.25cm}\\
\multicolumn{7}{|c|}{$y'(-1)=y(1)=0, \qquad p = p_{L} = 6-4\times(1/4)=5$}\\
\multicolumn{7}{|c|}{~}\vspace{-.25cm}\\
\hline
&&&&&&\vspace{-.25cm}\\
$N$ & $\delta\lambda_{5,N}$ & order & $\delta\lambda_{10,N}$ & order & $\delta\lambda_{20,N}$ & order \\
&&&&&&\vspace{-.25cm}\\
\hline
&&&&&&\vspace{-.25cm}\\
$49$  & $4.4416{\rm E}-06$ & $5.002$ & $5.5319{\rm E}-06$ & $4.993$ & $4.8732{\rm E}-06$ & $4.969$ \\ 
$99$  & $1.3859{\rm E}-07$ & $5.001$ & $1.7368{\rm E}-07$ & $4.999$ & $1.5556{\rm E}-07$ & $4.996$ \\ 
$199$ & $4.3295{\rm E}-09$ & $4.999$ & $5.4312{\rm E}-09$ & $4.997$ & $4.8743{\rm E}-09$ & $5.000$ \\ 
$399$ & $1.3534{\rm E}-10$ &   --    & $1.7002{\rm E}-10$ &   --    & $1.5234{\rm E}-10$ &   --   \\ 
\hline
\multicolumn{7}{|c|}{~}\vspace{-.25cm}\\
\multicolumn{7}{|c|}{$y(-1)=y'(1)=0, \qquad p = p_{R} = 6-4\times(3/4)=3$}\\
\multicolumn{7}{|c|}{~}\vspace{-.25cm}\\
\hline
&&&&&&\vspace{-.25cm}\\
$N$ & $\delta\lambda_{5,N}$ & order & $\delta\lambda_{10,N}$ & order & $\delta\lambda_{20,N}$ & order \\
&&&&&&\vspace{-.25cm}\\
\hline
&&&&&&\vspace{-.25cm}\\
$49$  & $2.1678{\rm E}-03$ & $3.002$ & $7.9981{\rm E}-03$ & $2.999$ & $1.2424{\rm E}-02$ & $2.985$ \\ 
$99$  & $2.7065{\rm E}-04$ & $3.000$ & $1.0005{\rm E}-03$ & $3.000$ & $1.5693{\rm E}-03$ & $2.999$ \\ 
$199$ & $3.3820{\rm E}-05$ & $3.000$ & $1.2505{\rm E}-04$ & $3.000$ & $1.9636{\rm E}-04$ & $3.000$ \\ 
$399$ & $4.2273{\rm E}-06$ &   --    & $1.5631{\rm E}-05$ &   --    & $2.4547{\rm E}-05$ &   --   \\ 
\hline
\multicolumn{7}{|c|}{~}\vspace{-.25cm}\\
\multicolumn{7}{|c|}{$y'(\pm 1)=0, \qquad p = p_{R} = 3$}\\
\multicolumn{7}{|c|}{~}\vspace{-.25cm}\\
\hline
&&&&&&\vspace{-.25cm}\\
$N$ & $\delta\lambda_{5,N}$ & order & $\delta\lambda_{10,N}$ & order & $\delta\lambda_{20,N}$ & order \\
&&&&&&\vspace{-.25cm}\\
\hline
&&&&&&\vspace{-.25cm}\\
$49$  & $1.7360{\rm E}-03$ & $3.001$ & $7.4699{\rm E}-03$ & $3.003$ & $1.2385{\rm E}-02$ & $2.989$ \\ 
$99$  & $2.1688{\rm E}-04$ & $3.000$ & $9.3168{\rm E}-04$ & $3.001$ & $1.5602{\rm E}-03$ & $2.999$ \\ 
$199$ & $2.7105{\rm E}-05$ & $3.000$ & $1.1637{\rm E}-04$ & $3.000$ & $1.9510{\rm E}-04$ & $3.000$ \\ 
$399$ & $3.3881{\rm E}-06$ &   --    & $1.4543{\rm E}-05$ &   --    & $2.4386{\rm E}-05$ &   --   \\ 
\hline              
\multicolumn{7}{|c|}{~}\vspace{-.25cm}\\
\multicolumn{7}{|c|}{$y'(\pm 1)= y(\pm 1), \qquad p = p_{R} = 3$}\\
\multicolumn{7}{|c|}{~}\vspace{-.25cm}\\
\hline
&&&&&&\vspace{-.25cm}\\
$N$ & $\delta\lambda_{5,N}$ & order & $\delta\lambda_{10,N}$ & order & $\delta\lambda_{20,N}$ & order \\
&&&&&&\vspace{-.25cm}\\
\hline
&&&&&&\vspace{-.25cm}\\
$49$  & $2.0295{\rm E}-03$ & $3.002$ & $8.1569{\rm E}-03$ & $3.004$ & $1.2697{\rm E}-02$ & $2.990$ \\ 
$99$  & $2.5338{\rm E}-04$ & $3.000$ & $1.0168{\rm E}-03$ & $3.001$ & $1.5983{\rm E}-03$ & $3.000$ \\ 
$199$ & $3.1663{\rm E}-05$ & $3.000$ & $1.2698{\rm E}-04$ & $3.000$ & $1.9983{\rm E}-04$ & $3.000$ \\ 
$399$ & $3.9576{\rm E}-06$ &   --    & $1.5869{\rm E}-05$ &   --    & $2.4977{\rm E}-05$ &   --   \\ 
\hline              
\end{tabular}
\end{center}
\end{small}
\end{table}

\begin{table}[t]
\caption{Order of convergence for problems with potential (\ref{q2}).} \label{tab2}
\begin{small}
\begin{center}
\begin{tabular}{|r|c c|c c|c c|}\hline
\multicolumn{7}{|c|}{~}\vspace{-.25cm}\\
\multicolumn{7}{|c|}{$y'(-1)=y(1)=0, \qquad p = p_{L} = 6-4\times(1/2) = 4$}\\
\multicolumn{7}{|c|}{~}\vspace{-.25cm}\\
\hline
&&&&&&\vspace{-.25cm}\\
$N$ & $\delta\lambda_{5,N}$ & order & $\delta\lambda_{10,N}$ & order & $\delta\lambda_{20,N}$ & order \\
&&&&&&\vspace{-.25cm}\\
\hline
&&&&&&\vspace{-.25cm}\\
$49$  & $6.1520{\rm E-}05$ & $3.998$ & $7.5495{\rm E-}05$ & $3.993$ & $6.9334{\rm E-}05$ & $3.974$ \\ 
$99$  & $3.8510{\rm E-}06$ & $3.999$ & $4.7406{\rm E-}06$ & $3.999$ & $4.4108{\rm E-}06$ & $3.997$ \\ 
$199$ & $2.4080{\rm E-}07$ & $4.000$ & $2.9654{\rm E-}07$ & $4.000$ & $2.7632{\rm E-}07$ & $4.000$ \\ 
$399$ & $1.5052{\rm E-}08$ &   --    & $1.8538{\rm E-}08$ &   --    & $1.7275{\rm E-}08$ &   --   \\ 
\hline
\multicolumn{7}{|c|}{~}\vspace{-.25cm}\\
\multicolumn{7}{|c|}{$y(-1)=y'(1)=0, \qquad p = p_{R} = 6-4\times(7/8)= 2.5$}\\
\multicolumn{7}{|c|}{~}\vspace{-.25cm}\\
\hline
&&&&&&\vspace{-.25cm}\\
$N$ & $\delta\lambda_{5,N}$ & order & $\delta\lambda_{10,N}$ & order & $\delta\lambda_{20,N}$ & order \\
&&&&&&\vspace{-.25cm}\\
\hline
&&&&&&\vspace{-.25cm}\\
$49$  & $6.9840{\rm E-}03$ & $2.503$ & $2.4576{\rm E-}02$ & $2.502$ & $4.1088{\rm E-}02$ & $2.491$ \\ 
$99$  & $1.2317{\rm E-}03$ & $2.501$ & $4.3385{\rm E-}03$ & $2.501$ & $7.3081{\rm E-}03$ & $2.500$ \\ 
$199$ & $2.1760{\rm E-}04$ & $2.500$ & $7.6648{\rm E-}04$ & $2.500$ & $1.2921{\rm E-}03$ & $2.500$ \\ 
$399$ & $3.8460{\rm E-}05$ &   --    & $1.3547{\rm E-}04$ &   --    & $2.2839{\rm E-}04$ &   --   \\ 
\hline
\multicolumn{7}{|c|}{~}\vspace{-.25cm}\\
\multicolumn{7}{|c|}{$y'(\pm 1)=0, \qquad p = p_{R} =2.5$}\\
\multicolumn{7}{|c|}{~}\vspace{-.25cm}\\
\hline
&&&&&&\vspace{-.25cm}\\
$N$ & $\delta\lambda_{5,N}$ & order & $\delta\lambda_{10,N}$ & order & $\delta\lambda_{20,N}$ & order \\
&&&&&&\vspace{-.25cm}\\
\hline
&&&&&&\vspace{-.25cm}\\
$49$  & $5.6379{\rm E-}03$ & $2.507$ & $2.3091{\rm E-}02$ & $2.510$ & $4.0904{\rm E-}02$ & $2.497$ \\ 
$99$  & $9.9209{\rm E-}04$ & $2.502$ & $4.0545{\rm E-}03$ & $2.503$ & $7.2462{\rm E-}03$ & $2.501$ \\ 
$199$ & $1.7508{\rm E-}04$ & $2.501$ & $7.1512{\rm E-}04$ & $2.501$ & $1.2797{\rm E-}03$ & $2.501$ \\ 
$399$ & $3.0932{\rm E-}05$ &   --    & $1.2633{\rm E-}04$ &   --    & $2.2612{\rm E-}04$ &   --   \\ 
\hline              
\multicolumn{7}{|c|}{~}\vspace{-.25cm}\\
\multicolumn{7}{|c|}{$y'(\pm 1)= y(\pm 1), \qquad p = p_{R} =2.5$}\\
\multicolumn{7}{|c|}{~}\vspace{-.25cm}\\
\hline
&&&&&&\vspace{-.25cm}\\
$N$ & $\delta\lambda_{5,N}$ & order & $\delta\lambda_{10,N}$ & order & $\delta\lambda_{20,N}$ & order \\
&&&&&&\vspace{-.25cm}\\
\hline
&&&&&&\vspace{-.25cm}\\
$49$  & $6.4759{\rm E-}03$ & $2.506$ & $2.5025{\rm E-}02$ & $2.510$ & $4.2056{\rm E-}02$ & $2.498$ \\ 
$99$  & $1.1403{\rm E-}03$ & $2.502$ & $4.3929{\rm E-}03$ & $2.503$ & $7.4452{\rm E-}03$ & $2.502$ \\ 
$199$ & $2.0131{\rm E-}04$ & $2.501$ & $7.7477{\rm E-}04$ & $2.501$ & $1.3146{\rm E-}03$ & $2.501$ \\ 
$399$ & $3.5571{\rm E-}05$ &   --    & $1.3686{\rm E-}04$ &   --    & $2.3228{\rm E-}04$ &   --   \\ 
\hline              
\end{tabular}
\end{center}
\end{small}
\end{table}

\begin{figure}[t]
\begin{center}
 \includegraphics[width=.95\textwidth, height = .7\textwidth]{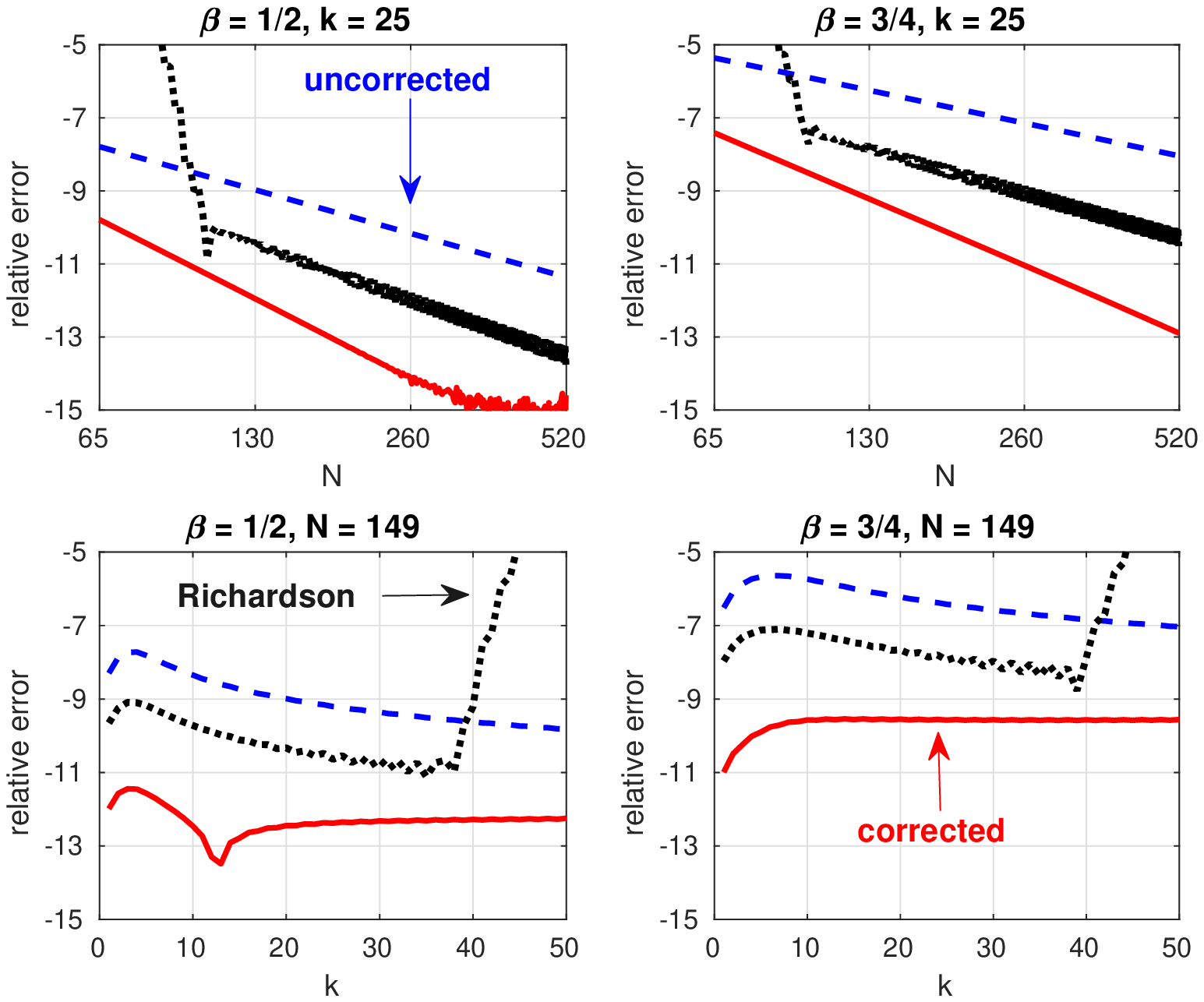}
\caption{Relative errors for the symmetric problems with potential (\ref{q3}) subject
to Neumann-Neumann BCs.}
 \label{fig2}
\end{center}
\end{figure}

\begin{figure}[t]
\begin{center}
 \includegraphics[width=.95\textwidth, height = .7\textwidth]{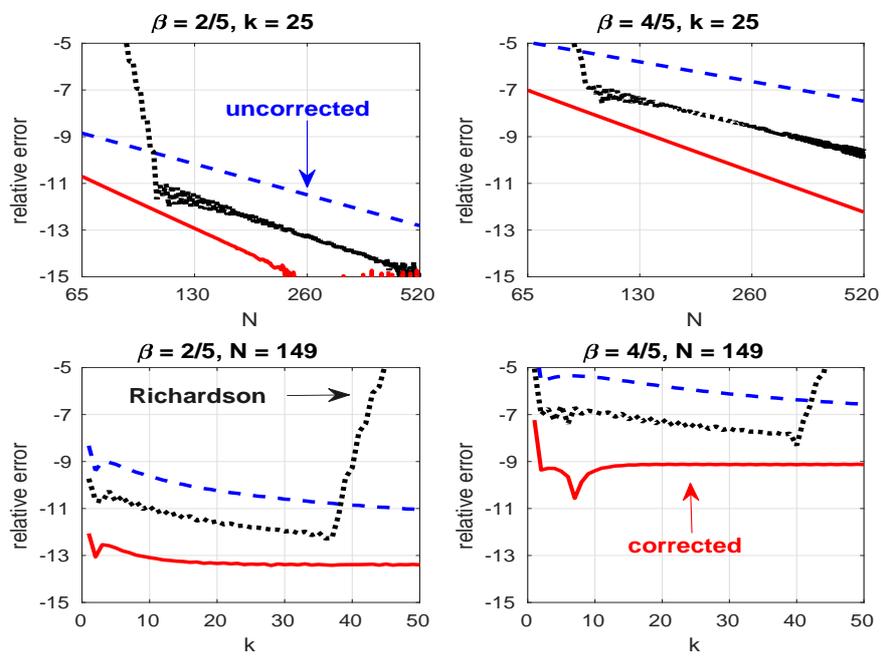}
\caption{Relative errors for problems with potentials (\ref{q4})
subject to $y(-1)-y'(-1)=y'(1)=0.$}
 \label{fig3}
\end{center}
\end{figure}

In the following final tests, we compare the performances of our correction technique with 
those of the classical Richardson extrapolation given by
$$\rho_k^{(N)} \equiv \frac{2^p\,\lambda_k^{(N)} - \lambda_k^{(N-1)/2}}{2^p-1} \approx \lambda_k, \qquad N \mbox{ odd},$$
where $p$ is specified in the convergence theorem.  In particular, we considered symmetric problems with  potentials
\begin{equation}\label{q3}
q(x) = \frac{10}{(2-x^2)(1-x^2)^\be}, \quad \be =\frac{1}2,\frac{3}4, 
\end{equation}
and problems with the following not symmetric $q$'s
\begin{equation}\label{q4}
q(x) = 5\,\frac{\cosh(x)\,(1+x)^{\frac{1}{5}}+2\log({\frac{3}{2}}+x)(1-x)^{\frac{1}{5}}+4(1-x^2)}{(1-x^2)^\be},\,\,\, \be = \frac{2}{5},\frac{4}5.
\end{equation}
The results obtained for some BCs have been reported in Figures~\ref{fig2} and \ref{fig3} (the errors corresponding 
to $\rho_k^{(N)}$ are depicted in dotted lines). As one can see, the Richardson extrapolation requires $N$ much larger than 
$k$ to improve the eigenvalue approximation, say $N$ not smaller than $4k+1.$  If this is not the case then it may deteriorates 
drastically the accuracy of $\lambda_k^{(N)}.$ Furthermore, the improvement that we get 
with our low-cost method is undeniably larger than that obtained with Richardson.

\end{document}